\documentclass[a4paper,12pt]{amsart}
\usepackage{amsfonts}
\usepackage{amssymb}
\usepackage{ifthen}
\usepackage{graphicx}
\nonstopmode \numberwithin{equation}{section}
\usepackage[usenames]{color}
\newtheorem{thm}{Theorem}[section]
\newtheorem{lem}{Lemma}[section]
\newtheorem{cor}{Corollary}[section]
\newtheorem{prop}{Proposition}[section]

\newtheorem{cl}{Claim}[section]
\newtheorem{ca}{Case}[section]
\newtheorem{sca}{Subcase}[section]
\newtheorem{scl}{Subclaim}[section]
\newtheorem{conj}[equation]{Conjecture}

\theoremstyle{definition}
\newtheorem{defn}{Definition}[section]
\newtheorem{op}[equation]{Open Problem}
\newtheorem{ques}[equation]{Question}
\newtheorem{rem}{Remark}[section]
\newtheorem{exam}[equation]{Example}

\newcounter {own}
\def\theown {\thesection       .\arabic{own}}

\newenvironment{pf}[1][]{%
 \vskip 3mm
 \noindent
 \ifthenelse{\equal{#1}{}}%
  {{\slshape Proof. }}%
  {{\slshape #1.} }%
 }%
{\qed\bigskip}

\newcounter{alphabet}
\newcounter{tmp}

\makeatletter
\newcommand{\Ref}[1]{\@ifundefined{r@#1}{}{\setcounter{tmp}{\ref{#1}}\Alph{tmp}}}
\makeatother
\newcommand{\IR}{{\mathbb R}}

\newcommand{\diam}{{\operatorname{diam}}}

\def\be{\begin{equation}}
\def\ee{\end{equation}}

\newcommand{\ben}{\begin{enumerate}}
\newcommand{\een}{\end{enumerate}}

\newcommand{\blem}{\begin{lem}}
\newcommand{\elem}{\end{lem}}
\newcommand{\bthm}{\begin{thm}}
\newcommand{\ethm}{\end{thm}}
\newcommand{\bcor}{\begin{cor}}
\newcommand{\ecor}{\end{cor}}
\newcommand{\beg}{\begin{exam}}
\newcommand{\eeg}{\end{exam}}
\newcommand{\begs}{\begin{examples}}
\newcommand{\eegs}{\end{examples}}
\newcommand{\bdefe}{\begin{defn}}
\newcommand{\edefe}{\end{defn}}
\newcommand{\bprob}{\begin{prob}}
\newcommand{\eprob}{\end{prob}}
\newcommand{\bques}{\begin{ques}}
\newcommand{\eques}{\end{ques}}
\newcommand{\bei}{\begin{itemize}}
\newcommand{\eei}{\end{itemize}}
\newcommand{\bcon}{\begin{conj}}
\newcommand{\econ}{\end{conj}}
\newcommand{\bop}{\begin{op}}
\newcommand{\eop}{\end{op}}

\newcommand{\bca}{\begin{ca}}
\newcommand{\eca}{\end{ca}}
\newcommand{\bsca}{\begin{sca}}
\newcommand{\esca}{\end{sca}}

\newcommand{\bcl}{\begin{cl}}
\newcommand{\ecl}{\end{cl}}

\newcommand{\bscl}{\begin{scl}}
\newcommand{\escl}{\end{scl}}

\newcommand{\bcons}{\begin{conjs}}
\newcommand{\econs}{\end{conjs}}

\newcommand{\bprop}{\begin{prop}}
\newcommand{\eprop}{\end{prop}}

\newcommand{\br}{\begin{rem}}
\newcommand{\er}{\end{rem}}
\newcommand{\brs}{\begin{rems}}
\newcommand{\ers}{\end{rems}}
\newcommand{\bo}{\begin{obser}}
\newcommand{\eo}{\end{obser}}
\newcommand{\bos}{\begin{obsers}}
\newcommand{\eos}{\end{obsers}}
\newcommand{\bpf}{\begin{pf}}
\newcommand{\epf}{\end{pf}}
\newcommand{\ba}{\begin{array}}
\newcommand{\ea}{\end{array}}
\newcommand{\beq}{\begin{eqnarray}}
\newcommand{\beqq}{\begin{eqnarray*}}
\newcommand{\eeq}{\end{eqnarray}}
\newcommand{\eeqq}{\end{eqnarray*}}

\newcounter{minutes}
\divide\time by 60
\newcounter{hours}
\multiply\time by 60 \addtocounter{minutes}{-\time}
\begin{document}

\bibliographystyle{amsplain}
\title{Uniform perfectness for quasi-metric spaces}

%%%%%%%% BEGIN TIMESTAMP
\def\thefootnote{}
\footnotetext{ \texttt{\tiny File:~\jobname .tex,
          printed: \number\year-\number\month-\number\day,
          \thehours.\ifnum\theminutes<10{0}\fi\theminutes}
} \makeatletter\def\thefootnote{\@arabic\c@footnote}\makeatother
%%%%%%%% END TIMESTAMP

\author{Qingshan Zhou}
\address{Qingshan Zhou,  Department of Mathematics, Shantou University, Shantou 515063, People's Republic
of China} \email{q476308142@qq.com}

\author{Yaxiang  Li${}^{\mathbf{*}}$}
\address{Yaxiang Li,  Department of Mathematics, Hunan First Normal University, Changsha,
Hunan 410205, People's Republic
of China}\address{Hunan Provincial Key Laboratory of Mathematical Modeling and Analysis in Engineering, Changsha University of Science and Technology, Changsha 410114, Hunan, People's Republic
of China} \email{yaxiangli@163.com}

\author{Ailing Xiao}
\address{Ailing Xiao, School of Finance, Guangdong University of Foreign Studies, Guangzhou 510006, People's Republic
of China} \email{alxiao181@126.com}

\date{}
 \keywords{Uniform perfectness, quasi-metric space, quasim\"obius map, quasisymmetric map, Gromov hyperbolic space.\\
${}^{\mathbf{*}}$ Corresponding author}

\begin{abstract}
The aim of this paper is to investigate the equivalence conditions for uniform perfectness of quasi-metric spaces. We also obtain the invariant property of uniform perfectness under quasim\"obius maps in quasi-metric spaces. In the end, two applications are given.
\end{abstract}

\thanks{The research was partly supported by NNSF of
China (Nos. 11601529,  11671127,  11571216) and Hunan Provincial Key Laboratory of Mathematical Modeling and Analysis in Engineering  (No.2017TP1017).}

\maketitle{} \pagestyle{myheadings} \markboth{}{}

\section{Introduction and main results}\label{sec-1}
In this paper, we mainly study several equivalence conditions for uniform perfectness of quasi-metric spaces. We start with the definition of quasi-metric spaces.
\bdefe For a given set $Z$ and a constant $K\geq 1$,
\ben
\item
a function $\rho:$ $Z\times Z\to [0, +\infty)$ is said to be {\it $K$-quasi-metric}
\ben
\item
if
for all $x$ and $y$ in $Z$, $\rho(x,y)\geq 0$, and $\rho(x,y)=0$ if and only if $x=y$;
\item
 $\rho(x,y)=\rho(y,x)$ for all $x,$ $y\in Z$;

\item
 $\rho(x,z)\leq K(\rho(x,y)\vee\rho(y,z))$ for all $x,$ $y,$ $z\in Z$.
 \een
\item
the pair $(Z, \rho)$ is said to be a {\it $K$-quasi-metric space} if the function $\rho:$ $Z\times Z\to [0, +\infty)$ is $K$-quasi-metric with $K\geq 1$. In the following, we always say that $K$ is the quasi-metric coefficient of $(Z, \rho)$.
\een
Here and hereafter, we use the notations: $r\vee s$ and $r\wedge s$ for numbers $r,$ $s$ in $\mathbb{R}$, where $$r\vee s=\max\{r,s\}\;\;\mbox{ and}\;\;r\wedge s=\min\{r,s\}.$$
\edefe
Clearly, if $(Z, \rho)$ is $K_1$-quasi-metric, it must be $K_2$-quasi-metric for any $K_2\geq K_1$. Hence, in the following, the quasi-metric coefficients of all quasi-metric spaces are always assumed to be $K$ with $K> 1$. For more properties and examples concerning quasi-metric spaces, see \cite{Ass, BuSc, CW, FR37, Hei, Ibra, Ibra02, JJ, MS} etc.

Every quasi-metric $\rho$ defines a uniform structure on $Z$. The balls
$\mathbb{B}(x, r) =\{y\in Z:\; \rho(x,y)<r\}$ $(r > 0)$ form a basis of neighborhoods of $x$ for the topology
induced by the uniformity on $Z$. We shall refer to this
topology as the $\rho$-topology of $Z$ (cf. \cite{MS}).

%{\color{red} Note that the quasi-metric spaces generally have no fixed topology, and even if it does, we do not know whether the above quasi-metric balls are open or %not. In this paper, we shall refer to the topology of $(Z,\rho)$ by using the quasi-metric balls, $\mathbb{B}_{\rho}(x,r)$, $r>0$, which actually form a basis of %neighborhoods of $x$ for the topology induced by the uniformity on $Z$; see \cite{MS}. }

Quasi-metric spaces are not quite as well behaved as metric spaces. For instance, the unit disk in the plane with a power (which is less than one) of the Euclidean metric is clearly a quasi-metric space and one can check that there are no rectifiable curves in this space except the constant curves.
Also, ``open balls" in quasi-metric spaces may be not open in topology associated with quasi-metric.
The following useful result on the relationship between quasi-metric spaces and metric spaces follows from \cite[Proposition 2.2.5]{BuSc}.
\blem\label{lem4}
Let $(X,\rho)$ be a $K$-quasi-metric space. For a constant $0<\varepsilon\leq 1$, if $K^{\varepsilon}\leq 2$, then there is a metric $d_{\varepsilon}$ on $X$ such that $$\frac{1}{4}\rho^{\varepsilon}(z_1,z_2)\leq d_{\varepsilon}(z_1,z_2)\leq \rho^{\varepsilon}(z_1,z_2) $$ for all $z_1,z_2\in X$.
\elem

To state our results, we need the definition of uniform perfectness in quasi-metric spaces.

\bdefe
A quasi-metric space $(Z,\rho)$ is called {\it uniformly perfect} if there is a constant $\mu\in (0,1)$ such that for each $x\in Z$ and every $r>0$, the set $\mathbb{B}(x,r)\setminus \mathbb{B}(x, \mu r)\not=\emptyset$ provided that $Z\setminus \mathbb{B}(x,r)\not=\emptyset$.
\edefe

Uniform perfectness is a weaker condition than connectedness. Connected spaces are uniformly perfect, and those with isolated points are not. Many disconnected fractals such as the Cantor ternary set, Julia sets and the limit set of a nonelementary, finitely generated Kleinian group of $\overline{\IR}^n$ are uniformly perfect \cite{JV}. In particular, uniform perfectness has provided a useful tool in modern research of geometric functional theory, harmonic analysis and asymptotic geometry; see \cite{BuSc, CW, Heer, Hei, JJ, Ve, WZ, ZLL}.

It is worth mentioning that Buyalo and Schroeder established the quasisymmetric and quasim\"{o}bius extension theorems for visual geodesic hyperbolic spaces which possess uniformly perfect boundaries at infinity \cite[Chapter $7$]{BuSc}. In \cite{WZ}, the first author and Wang found  several conditions under which a weakly quasim\"obius map is quasim\"obius in uniformly perfect quasi-metric spaces. The authors in \cite{ZLL} investigated the invariance of doubling property under sphericalization and flattening transformations in uniformly perfect spaces. Recently, Vellis \cite{Ve} proved the classical quasisymmetric Schoenflies theorem for planar uniform domains with uniformly prefect boundaries.

As the first goal of this paper, we shall discuss the relationship among uniform perfectness, homogeneous density, $\sigma$-density etc in the setting of quasi-metric spaces. We show that all these conditions are equivalent. (Note that other notions appearing in the following results will be introduced in the body of this paper.)

 \begin{thm}\label{main-thm1} Suppose $(Z, \rho)$ is a quasi-metric space. Then the following are quantitatively equivalent.
 \ben
\item\label{fri-1}
 $(Z, \rho)$ is $\mu$-uniformly perfect;
\item\label{fri-2}
$(Z, \rho)$ is $(\lambda_1,\lambda_2)$-HD;
\item\label{fri-3}
$(Z, \rho)$ is $\sigma$-dense;
\item\label{fri-4}
There are numbers $\mu_1$ and $\mu_2$ such that $0< \mu_1\leq \mu_2<1$ and for any triple $(a,c,d)$ in $Z$, there is a point $x\in Z$ satisfying $\mu_1\leq r(a,x,c,d)\leq \mu_2$.
\een
\end{thm}
Here, we make the following notational convention: Suppose $A$ denotes a condition with data $v$ and $B$ another condition with data $v'$. We say that $A$ implies $B$ quantitatively if $A$ implies $B$ so that $v'$ depends only on $v$. If $A$ and $B$ imply each other quantitatively, then we say that they are quantitatively equivalent.

\br
In metric spaces, the equivalence between \eqref{fri-1} and \eqref{fri-4} (resp. \eqref{fri-2} and \eqref{fri-3}, \eqref{fri-3} and \eqref{fri-4}) coincides with \cite[Lemma 11.7]{Hei} (resp. \cite[Lemma 3.1]{Ibra02}, \cite[Remark 3.3]{Ase})
\er

In \cite{Vai89}, V\"{a}is\"{a}l\"{a} introduced the term quasim\"{o}bius maps in metric spaces and obtained the close connections among quasim\"{o}bius maps, quasiconformal maps and quasisymmetric maps. See \cite{Ase, BK, BuSc, Hei, TV, Vai-4} for more background materials in this line. Further, in \cite{TV}, Tukia and V\"{a}is\"{a}l\"{a} proved that every quasisymmetric map in uniformly prefect spaces is power quasisymmetric. As the second goal of this paper, we shall investigate the relationships of the uniform perfectness with (power) quasisymmetric maps and (power) quasim\"obius maps, respectively. One of our results reads as follows.

\begin{thm}\label{lem3}
Suppose $f:$ $(Z_1, \rho_1)\to (Z_2, \rho_2)$ is a quasim\"obius map between quasi-metric spaces $(Z_1, \rho_1)$ and $(Z_2, \rho_2)$. Then $(Z_1, \rho_1)$ is uniformly perfect if and only if  $(Z_2, \rho_2)$ is uniformly perfect, quantitatively.
\end{thm}

We remark that Theorem \ref{lem3} coincides with \cite[Corollary 4.6]{JV} in the setting of $\IR^n$. With the help of Theorem \ref{lem3}, we arrive at the following equivalent conditions for uniform perfectness concerning quasisymmetric and quasim\"obius maps in quasi-metric spaces.

\begin{thm}\label{main-thm0} Suppose $(Z, \rho)$ is a quasi-metric space. If $(Z,\rho)$ has no isolated points, then\ben
\item\label{main-1}
$(Z,\rho)$ is uniformly perfect if and only if every quasisymmetric map of $(Z,\rho)$ to a quasi-metric space is power quasisymmetric, quantitatively;
\item\label{main-2}
$(Z,\rho)$ is uniformly perfect if and only if every quasim\"obius map of $(Z,\rho)$ to a quasi-metric space is power quasim\"obius, quantitatively.
\een
\end{thm}

\br $(1)$ In metric spaces, Theorem \ref{main-thm0}\eqref{main-1} is the same as \cite[Theorems 4.13 and 6.20]{TroVai}.
$(2)$
Aseev and Trotsenko proved that, in $\IR^n$, if $(Z,\rho)$ is $\sigma$-dense, then every quasim\"obius map of $(Z,\rho)$ is power quasim\"obius by applying the conformal
moduli of families of curves (see \cite[Theorem 4.1]{Ase}). One finds that the method of proof in \cite{Ase} is no longer valid in metric spaces. So, even in $\IR^n$, the method of proof of Theorem \ref{main-thm0}\eqref{main-2} is new.
\er

Recently, Meyer studied the relationship between Gromov hyperbolic spaces and their boundaries at infinity. He proved the invariant property of the uniform perfectness with respect to the inversions in quasi-metric spaces (see \cite[Theorem 7.1]{Me}). This result is one of the main results in \cite{Me}, whose proof is lengthy.
As an application of Theorem \ref{lem3}, we shall give a different proof to \cite[Theorem 7.1]{Me} (see Theorem \ref{app-1} below). Also, we shall discuss the uniform perfectness of a complete quasi-metric space and the corresponding boundary of its hyperbolic approximation by applying Theorem \ref{lem3} (see Theorem \ref{app-2} below).

We remark that after this work was first finished in 2016, there were several applications based on that version. Based on Theorems \ref{main-thm1} and \ref{lem3}, the first author and  Wang studied the relations between quasim\"obius mappings and weakly quasim\"obius mappings \cite{WZ}. And then by using the main results in \cite{WZ}, Vellis studied the extension properties of planer uniform domains \cite{Ve},  and recently Aseev studied the BAD class of multivalued mappings of Ptolemaic M\"obius structures in the sense of Buyalo with controlled distortion of generalized angles \cite{Ase2018}.

The organization of this paper is as follows. In the second section, we shall introduce some necessary concepts and discuss the condition in quasi-metric spaces under which quasim\"obius maps and quasisymmetric maps are the same. In the third section, some concepts will be introduced and the equivalence of uniform perfectness with homogeneous density, $\sigma$-density etc will be proved. The invariant property of uniform perfectness  with respect to quasim\"obius maps will be shown in the forth section, and in the fifth section, relationships among (power) quasisymmetric maps, (power) quasim\"obius maps and uniform perfectness will be established. In the last section, some applications of Theorem \ref{lem3} will be given.

%%%%%%%%%%%%%%%%%%%%%%%%%%%%%%%
%%%%%%%%%%%%%%%%%%%%%%%%%%%%%%%
\section{Quasim\"obius maps and quasisymmetric maps in quasi-metric spaces}\label{sec-2}
%%%%%%%%%%%%%%%%%%%%%%%%%%%%%%%
%%%%%%%%%%%%%%%%%%%%%%%%%%%%%%%

In this section, we shall introduce some necessary notations and concepts, and prove several basic results. The main result in this section is Theorem \ref{w-04}, which concerns the condition under which power quasisymmetric maps and power quasim\"obius maps are the same.

%%%%%%%%%%%%%%%%%%%%%%%%%%%%%%%
\subsection{Cross ratios}\ \ \
%%%%%%%%%%%%%%%%%%%%%%%%%%%%%%%

For four points $a,$ $b,$ $c,$ $d$ in a quasi-metric space $(Z,\rho)$, its {\it cross ratio} is defined by the number
$$r(a,b,c,d)=\frac{\rho(a,c)\rho(b,d)}{\rho(a,b)\rho(c,d)}.$$
Then we have
\bprop\label{mon-7}
\noindent $(1)$ For $a,$ $b,$ $c$ and $d$ in $(Z,\rho)$, $$r(a,b,c,d)=\frac{1}{r(b, d, a, c)};$$

\noindent $(2)$ For $a,$ $b,$ $c,$ $d$ and $z$ in $(Z,\rho)$, $$r(a,b,c,d)=r(a,b,z,d)r(a,z,c,d).$$
\eprop

In \cite{BK}, Bonk and  Kleiner introduced the following useful notation:
$$\langle a,b,c,d \rangle=\frac{\rho(a,c)\wedge\rho(b,d)}{\rho(a,b)\wedge\rho(c,d)}.$$

Bonk and Kleiner established a relation between $r(a,b,c,d)$ and $\langle a,b,c,d \rangle$ in the setting of metric spaces (see \cite[Lemma 3.3]{BK}).
The following result shows that this relation also holds in quasi-metric spaces.

\begin{lem}\label{w-01} For any $a,$ $b,$ $c,$ $d$ in $(Z,\rho)$, we have\ben
\item\label{mon-3}
$\frac{1}{\theta_K(1/r(a,b,c,d))} \leq \langle a,b,c,d \rangle \leq \theta_K(r(a,b,c,d));$
\item\label{mon-4}
$\theta_K^{-1}(\langle a,b,c,d \rangle)\leq r(a,b,c,d)  \leq
\frac{1}{\theta_K^{-1}(1/\langle a,b,c,d \rangle)} ,$
\een
where $\theta_K(t)=K^2(t\vee \sqrt{t})$. (Here, we recall that $K$ denotes the coefficient of the quasi-metric space $(Z, \rho)$.)
\end{lem}

\bpf Obviously, we only need to prove \eqref{mon-3} in the lemma. For the proof, we let
$$\langle a,b,c,d \rangle =s\;\;\mbox{and} \;\;r(a,b,c,d)=t.$$ Without loss of generality, we may assume $\rho(a,c)\leq \rho(b,d)$. Then we have
\begin{eqnarray*} \rho(a,b) &\leq& K(\rho(a,c)\vee\rho(c,b))
\leq K^2\big(\rho(a,c)\vee\rho(c,d)\vee\rho(d,b)\big)
\\ \nonumber&=& K^2\big(\rho(c,d)\vee\rho(d,b)\big),\end{eqnarray*}
and similarly, $$\rho(c,d)\leq K^2\big(\rho(a,b))\vee\rho(d,b)\big).$$ The combination of these two estimates leads to
$$ \rho(a,b)\vee\rho(c,d) \leq K^2\big((\rho(a,b)\wedge\rho(c,d)) \vee\rho(b,d)\big)\leq K^2(1\vee\frac{1}{s})\rho(d,b),
$$
and so we get
$$t=r(a,b,c,d)=\frac{\rho(a,c)\rho(b,d)}{\big(\rho(a,b)\wedge\rho(c,d)\big)\big(\rho(a,b)\vee\rho(c,d)\big)}\geq \frac{s}{K^2(1\vee\frac{1}{s})},$$
which implies $$s\leq \theta_K(t)=K^2(t\vee \sqrt{t}).$$ Hence the right side inequality in \eqref{mon-3} holds.

By Proposition \ref{mon-7}, we see that the left side inequality in \eqref{mon-3} easily follows from the right side one, and so the proof of the lemma is complete.
\epf

%%%%%%%%%%%%%%%%%%%%%%%%%%%%%%%
\subsection{Quasisymmetric maps  and quasim\"obius maps}
%%%%%%%%%%%%%%%%%%%%%%%%%%%%%%%

\bdefe \label{def1-0} Suppose
$\eta$ and $\theta$ are homeomorphisms from $[0, \infty)$ to $[0, \infty)$.  A homeomorphism $f:$ $(Z_1,\rho_1)\to (Z_2,\rho_2)$ is said to be
\ben
\item\label{sun-1}
\ben
\item\label{tues-1}
{\it
$\eta$-quasisymmetric} if $\rho_1(x,a)\leq t \rho_1(x,b)$ implies $$\rho_2(x',a')\leq \eta(t) \rho_2(x',b')$$ for all $a,$ $b,$ $x$ in $(Z_1,\rho_1)$ and $t\geq 0$, where primes mean the images of points under $f$, for example, $x'=f(x)$ etc;
\item
{\it power quasisymmetric} if it is $\eta$-quasisymmetric, where $\eta$ has the form
$$\eta(t)=M_1(t^{1/\alpha}\vee t^{\alpha}) $$ for some constants $\alpha \geq 1$ and $M_1\geq 1$.
\een
\item\label{tues-2}
\ben
\item\label{sun-2}
{\it
$\theta$-quasim\"obius} if $r(a,b,c,d)\leq t$ implies $$ r(a',b',c',d')\leq \theta(t)$$
for all $a,$ $b,$ $c,$ $d$ in $(Z_1,\rho_1)$ and $t\geq 0$;
\item
{\it power quasim\"obius} if it is $\theta$-quasim\"obius, where $\theta$ has the form
$$\theta(t)=M_2(t^{1/\beta}\vee t^{\beta})$$
for some constants $\beta \geq 1$ and $M_2\geq 1$.
\een
\een
\edefe

As a direct consequence of Lemma \ref{w-01}, we have the following two results.

\begin{lem}\label{w-02}
Suppose $f:$ $(Z_1,\rho_1)\to (Z_2,\rho_2)$ is a homeomorphism between two quasi-metric spaces.
\ben
\item\label{mon-8}
If $f$ is $\eta$-quasisymmetric, then it is $\theta$-quasim\"{o}bius, where $\theta(t)=\frac{1}{\theta_{K}^{-1}(\frac{1}{\eta\circ\theta_K(t)})}$ and $\theta_{K}$ is from Lemma \ref{w-01}.
\item\label{thur-1}
If $f$ is a power quasisymmetric map with its control function $\eta(t)=M(t^{\alpha}\vee t^{\frac{1}{\alpha}})$, where $M\geq 1$ and $\alpha\geq 1$, then it is power quasim\"{o}bius with its control function $\theta(t)=M^2K^{4(1+\alpha)}(t^{2\alpha}\vee t^{\frac{1}{2\alpha}})$.
\een
\end{lem}

We remark that Lemma \ref{w-02}\eqref{mon-8} is a generalization of
\cite[Theorem 3.2]{Vai89} in the setting of quasi-metric spaces.

Next, we consider the converse of Lemma \ref{w-02}\eqref{thur-1} in the setting of bounded quasi-metric spaces. We shall discuss the converse of Lemma \ref{w-02}\eqref{mon-8} elsewhere. We start with the introduction of the following condition.

\bdefe
Suppose both $(Z_1,\rho_1)$ and $(Z_2,\rho_2)$ are bounded quasi-metric spaces. Let $\lambda\geq 1$ be a constant.
A homeomorphism $f:$ $(Z_1,\rho_1)\to (Z_2,\rho_2)$ is said to satisfy the {\it  $\lambda$-three-point condition} if there are points $z_1,$ $z_2,$ $z_3$ in $(Z_1, \rho_1)$ such that
$$\rho_1(z_i,z_j)\geq \frac{1}{\lambda}\diam(Z_1)\;\;\mbox{and}\;\; \rho_2(z'_i,z'_j)\geq \frac{1}{\lambda}\diam(Z_2)$$
for all $i\neq j\in\{1,2,3\}$, where ``$\diam$" means ``diameter".
\edefe

\begin{thm}\label{w-04}
Suppose that both $(Z_1,\rho_1)$ and $(Z_2,\rho_2)$ are bounded quasi-metric spaces and that $f:$ $(Z_1,\rho_1)\to (Z_2,\rho_2)$ satisfies the $\lambda$-three-point condition. Then
$f$ is power quasisymmetric if and only if it is power quasim\"{o}bius, quantitatively.
\end{thm}

\bpf
The necessity of the theorem obviously follows from Lemma \ref{w-02}\eqref{thur-1}. In the following, we prove the sufficiency. Let $f:$ $(Z_1,\rho_1)\to (Z_2,\rho_2)$ be a power quasim\"{o}bius map between two bounded quasi-metric spaces, which satisfies the $\lambda$-three-point condition for some constant $\lambda\geq 1$ and points $z_1,$ $z_2,$ $z_3\in Z_1$. We assume that the control function of $f$ is $$\theta(t)=M(t^{1/\beta}\vee t^{\beta})$$
for some constants $M\geq 1$ and $\beta \geq 1$.

To prove the power quasisymmetry of $f$, let $x,$ $a,$ $b$ be any three points in $(Z_1, \rho_1)$ with $\rho_1(x,a)=t \rho_1(x,b)$ with $t\geq 0$. Then we shall show that $$\rho_2(x',a')\leq \eta(t) \rho_2(x',b'),$$ where $\eta(t)=K^{3+6\beta}M(2\lambda)^{1+2\beta}(t^{1/(2\beta)}\vee t^{2\beta})$.

It follows from the $\lambda$-three-point condition that for any $w\in Z_1$,
 there are $i\not=j\in \{1, 2, 3\}$ such that
$$\rho_1(w,z_i)\wedge \rho_1(w,z_j)\geq \frac{\diam(Z_1)}{2K\lambda}.$$
Similarly, for any $u'\in Z_2$,
 there exist $m\not=n\in \{1, 2, 3\}$ such that
$$\rho_2(u',z_m')\wedge \rho_2(u',z_n')\geq \frac{\diam(Z_2)}{2K\lambda}.$$
Therefore, there must exist $z_i\in \{z_1,z_2,z_3\}$ such that $$\rho_1(a,z_i)\geq \frac{\diam(Z_1)}{2K\lambda}\;\;\mbox{and}\;\;\rho_2(b',z_i')\geq \frac{\diam(Z_2)}{2K\lambda}.$$ Thus
\beqq\label{w-06} \rho_1(a,z_i)\wedge \rho_1(x,b)\geq \frac{\rho_1(x,b)}{2K\lambda}\;\;\mbox{and}\;\;\rho_2(b',z_i')\wedge\rho_2(x',a')\geq \frac{\rho_2(x',a')}{2K\lambda},\eeqq
from which we deduce that
\be\label{w-07}\langle x,b,a,z_i\rangle\leq 2K\lambda \frac{\rho_1(x,a)}{\rho_1(x,b)}\;\;\mbox{and}\;\;  \langle x',b',a',z_i'\rangle\geq \frac{\rho_2(x',a')}{2K\lambda\rho_2(x',b')}.\ee

On the other hand, since $f$ is power quasim\"{o}bius with its control function $\theta$, we see from Lemma \ref{w-01} that
$$\langle x',b',a',z_i'\rangle \leq \theta'(\langle x,b,a,z_i\rangle),$$
where $\theta'(t)=\theta_K\circ \theta \Big(\frac{1}
{\theta_K^{-1}(1/ t)}\Big)$ and $\theta_K$
is from Lemma \ref{w-01}.
Then we deduce from (\ref{w-07}) that
$$\frac{\rho_2(x',a')}{\rho_2(x',b')}\leq 2K\lambda \theta'(\langle x,b,a,z_i\rangle)\leq 2K\lambda \theta'\Big(2K\lambda\frac{\rho_1(x,a)}{\rho_1(x,b)}\Big).$$
By taking $\eta(t)=K^{3+6\beta}M(2\lambda)^{1+2\beta}(t^{1/(2\beta)}\vee t^{2\beta})$, we see from elementary computations that
$$\frac{\rho_2(x',a')}{\rho_2(x',b')}\leq \eta\Big(\frac{\rho_1(x,a)}{\rho_1(x,b)}\Big).$$ Hence the proof of the theorem is complete.
\epf

\blem\label{sun-4} Suppose $f:$ $(Z_1,\rho_1)\to (Z_2,\rho_2)$ and $g:$ $(Z_2,\rho_2)\to (Z_3,\rho_3)$ are homeomorphisms. \ben
\item\label{sun-7}
If $f$ is $\theta_1$-quasim\"{o}bius and $g$ is $\theta_2$-quasim\"{o}bius, then $g\circ f$ is $\theta$-quasim\"{o}bius, where $\theta=\theta_{2}\circ \theta_1$;
\item\label{sun-5}
If $f$ is $\theta$-quasim\"{o}bius and $g$ is $\eta$-quasisymmetric, then $g\circ f$ is $\theta_1$-quasim\"{o}bius, where $\theta_1(t)=\frac{1}{\theta_{K}^{-1}\big(\frac{1}{\eta\circ\theta_K\circ \theta(t)}\big)}$;
\item\label{sun-6}
If $f$ is power quasim\"{o}bius and $g$ is power quasisymmetric (or power quasim\"obius), then $g\circ f$ is power quasim\"{o}bius, quantitatively.
\een
\elem

%%%%%%%%%%%%%%%%%%%%%%%%%%%%%%%
%%%%%%%%%%%%%%%%%%%%%%%%%%%%%%%
\section{uniform perfectness, homogeneous density and $\sigma$-density}\label{sec-3}
%%%%%%%%%%%%%%%%%%%%%%%%%%%%%%%
%%%%%%%%%%%%%%%%%%%%%%%%%%%%%%%

We start this section with several definitions, and then establish the invariant property of uniform perfectness with respect to quasisymmetric maps in quasi-metric spaces (Lemma \ref{lem2} below). Based on this result, Theorem \ref{main-thm1} will be proved.

%%%%%%%%%%%%%%%%%%%%%%%%%%%%%%%
\subsection{Homogeneous density and $\sigma$-density}\ \ \
%%%%%%%%%%%%%%%%%%%%%%%%%%%%%%%

\bdefe
Suppose $\{x_i\}_{i\in \mathbb{Z}}$ denotes a sequence of points in a quasi-metric space $(Z, \rho)$ with $a\not=x_i\not=b$.\ben
 \item[$(i)$]
 If $x_i\to a$ as $i\to-\infty$ and $x_i\to b$ as $i\to +\infty$, then $\{x_i\}$ is called a {\it chain} joining $a$ and $b$;
 further, if there is a constant $\sigma>1$ such that for all $i$, $$|\log r(a,x_i,x_{i+1},b)|\leq \log\sigma,$$ then $\{x_i\}$ is called a {\it $\sigma$-chain}.
 \item[$(ii)$]
 $(Z, \rho)$ is said to be {\it $\sigma$-dense} ($\sigma>1$) if any pair of points in $(Z, \rho)$ can be joined by a $\sigma$-chain.
 \een
\edefe

We remark that a $\sigma$-dense space does not contain any isolated point, and also, every $\sigma$-dense space must be $\sigma'$-dense for any $\sigma'\geq \sigma$.

\bdefe
A quasi-metric space $(Z, \rho)$ is said to be {\it homogeneously dense}, abbreviated $HD$, if there are constants $\lambda_1$ are $\lambda_2$ with $0<\lambda_1\leq \lambda_2<1$ such that
 for each pair of points $a,$ $b\in Z$, there is $x\in Z$ satisfying
$$\lambda_1\rho(a,b)\leq \rho(a,x)\leq\lambda_2\rho(a,b).$$
  To emphasize the parameters, we also say that $(Z, \rho)$ is $(\lambda_1,\lambda_2)$-$HD$.
\edefe

\begin{lem}\label{lem1}
$(1)$ If a quasi-metric space is $(\lambda_1,\lambda_2)$-$HD$, then it is $(\lambda_1^n,\lambda_2^n)$-$HD$ for any $n\in \mathbb{N}^+=\{1, 2, \ldots,\}$.

$(2)$ Suppose that both $(Z_1, \rho_1)$ and $(Z_2, \rho_2)$ are quasi-metric spaces and that $f:$ $(Z_1, \rho_1)\to  (Z_2, \rho_2)$ is $\eta$-quasisymmetric.
If $(Z_1, \rho_1)$ is $(\lambda_1,\lambda_2)$-$HD$, then $(Z_2, \rho_2)$ is $(\mu_1,\mu_2)$-$HD$, where both $\mu_1$ and $\mu_2$
depend only on $\lambda_1$, $\lambda_2$ and $\eta$.
\end{lem}

We remark that, in the setting of metric spaces, Lemma \ref{lem1} coincides with \cite[Lemma 3.9]{TV}. Also the proof of Lemma \ref{lem1} is similar to that of \cite[Lemma 3.9]{TV}. We omit it here.

%%%%%%%%%%%%%%%%%%%%%%%%%%%%%%%
\subsection{The invariant property of uniform perfectness with respect to quasisymmetric maps}
%%%%%%%%%%%%%%%%%%%%%%%%%%%%%%%

It is known that uniform perfectness is an invariant with respect to quasisymmetric maps in metric spaces (cf. \cite[Exercise 11.2]{Hei}). In the following, we
prove that this fact is still valid in quasi-metric spaces.

\begin{lem}\label{lem2} Let $f:$ $(Z_1,\rho_1)\to (Z_2,\rho_2)$ be $\eta$-quasisymmetric, where both $(Z_i,\rho_i)$ ($i=1$, $2$) are quasi-metric. Then $(Z_1, \rho_1)$ is uniformly perfect if and only if $(Z_2,\rho_2)$ is uniformly perfect, quantitatively.
\end{lem}
\bpf  Since the inverse of a quasisymmetric map is also quasisymmetric, to prove the lemma, it suffices to show that the uniform perfectness of $(Z_1, \rho_1)$ implies the uniform perfectness of $(Z_2, \rho_2)$. Now, we assume that $(Z_1, \rho_1)$ is $\mu$-uniformly perfect for some $\mu\in(0,1)$. Then we shall show that $(Z_2, \rho_2)$ is uniformly perfect too. To reach this goal, it suffices to find a constant $\mu'\in(0,1)$ such that for any $z'\in Z_2$ and $r>0$, if $Z_2\setminus \mathbb{B}(z',r)\not=\emptyset$, then there is a point $u'$ in $(Z_2, \rho_2)$ such that
$$\mu'r\leq \rho_2(z', u')<r.$$

By the assumption $Z_2\setminus \mathbb{B}(z',r)\not=\emptyset$, we see that there is a point $u_0'\in Z_2$ such that \be\label{mon-12}\rho_2(z',u'_0)\geq r.\ee

Choose $0<\alpha<1$ small enough such that $\eta(\alpha)< 1$. Then there exists an integer $k$ such that
\be\label{tues-4}\eta(\alpha)^{k}\rho_2(z',u_0') <r \leq \eta(\alpha)^{k-1}\rho_2(z',u_0').\ee

Since $(Z_1, \rho_1)$ is $\mu$-uniformly perfect and $u_0\in Z_1\setminus \mathbb{B}(z, \alpha\rho_1(z,u_0))$, we see that $\mathbb{B}(z, \alpha\rho_1(z,u_0))\setminus \mathbb{B}(z, \mu\alpha\rho_1(z,u_0))\not=\emptyset$. So there is a point $u_1\in Z_1$ such that
$$\mu\alpha\rho_1(z,u_0))\leq \rho_1(z,u_1)< \alpha\rho_1(z,u_0). $$
Hence
\be\label{tues-3}\mu'\rho_2(z',u'_0)\leq \rho_2(z',u_1')< \eta(\alpha)\rho_2(z',u'_0),\ee
where $\mu'=\frac{1}{\eta(\frac{1}{\mu\alpha})}$.

If $\rho_2(z',u_1')<r$, then \eqref{mon-12} and \eqref{tues-3} lead to
$$\mu'r\leq \rho_2(z', u_1')<r.$$
At present, we can take $u'=u_1'$.

Now, we consider the case: \be\label{mon-13}\rho_2(z',u_1')\geq r.\ee Since $(Z_1, \rho_1)$ is $\mu$-uniformly perfect and $u_1\in Z_1\setminus \mathbb{B}(z, \alpha\rho_1(z,u_1))$, we see that $\mathbb{B}(z, \alpha\rho_1(z,u_1))\setminus \mathbb{B}(z, \mu\alpha\rho_1(z,u_1))\not=\emptyset$. So there is a point $u_2\in Z_1$ such that
$$\mu\alpha\rho_1(z,u_1))\leq \rho_1(z,u_2)< \alpha\rho_1(z,u_1). $$
Hence $$\mu'\rho_2(z',u_1')\leq \rho_2(z',u_2')< \eta(\alpha)\rho_2(z',u_1')< \eta(\alpha)^2\rho_2(z',u'_0).$$

If $\rho_2(z',u_2')<r$, then \eqref{mon-13} leads to $$\mu'r\leq \rho_2(z', u_2')<r.$$
Hence, we can take $u'=u_2'$.

Next, we consider the case: $$\rho_2(z',u_2')\geq r.$$ By repeating this procedure, we can reach the following conclusion: There is $u_k'\in Z_2$ such that \ben
\item
For any $i\in \{1, \ldots, k-1\}$, $\rho_2(z',u_i')\geq r$;
\item
$\mu'\rho_2(z',u_{k-1}')\leq \rho_2(z', u_k')<\eta(\alpha)^k\rho_2(z',u'_0).$
\een
Then \eqref{tues-4} guarantees that
$$\mu'r\leq \rho_2(z', u_k')<r.$$
By taking $u'=u_k'$, we finish the proof.
\epf

%%%%%%%%%%%%%%%%%%%%%%%%%%%%%%%
\subsection{The proof of Theorem \ref{main-thm1}}
%%%%%%%%%%%%%%%%%%%%%%%%%%%%%%%

By applying Lemmas \ref{lem4}, \ref{lem1} and \ref{lem2}, together with Lemma \ref{lem-3} below, we see that the equivalence between $(1)$ and $(2)$ (resp. between $(2)$ and $(3)$)
easily follows from \cite[Lemma 11.7]{Hei} (resp. \cite[Lemma 3.1]{Ibra02}).
Hence, to finish the proof, it remains to show the equivalence between $(3)$ and $(4)$, whose proof is as follows.
\medskip

$(3)\Rightarrow(4)$ $\;$ Assume that $(Z,\rho)$ is $\sigma$-dense. Let $a,$ $c,$ $d$ be three distinct points in $(Z,\rho)$. Then there is a $\sigma$-chain $\{x_i\}_{i\in \mathbb{Z}}$ in $(Z,\rho)$ joining $a$ and $d$ such that
\beq\label{lem3.4-eq-1}\frac{1}{\sigma}\leq r(a, x_i,x_{i+1},d)\leq\sigma.\eeq
To prove this implication, it suffices to show that there is an integer $k$ such that
\be\label{thu-1}\frac{1}{2\sigma^2}\leq r(a, x_{k-1},c,d)\leq \frac{1}{2\sigma}.\ee

For the proof, we let
$$k=\inf\{i\in \mathbb{Z}:\; r(a, x_{i},c,d)<\frac{1}{2\sigma^2}\}.$$ Since $\lim_{i\to +\infty}r(a, x_{i},c,d)=0$ and $\lim_{i\to -\infty}r(a, x_{i},c,d)=+\infty$, we see that $k$ is finite, and so
$$r(a, x_{k},c,d)<\frac{1}{2\sigma^2}\;\;\mbox{and}\;\;r(a, x_{k-1},c,d)\geq \frac{1}{2\sigma^2}.$$
Then \eqref{lem3.4-eq-1} implies
$$r(a,x_{k-1},c,d)=r(a,x_k,c,d)r(a,x_{k-1},x_k,d)<\frac{1}{2\sigma}.$$
Hence \eqref{thu-1} is true, and thus the implication from $(3)$ to $(4)$ is proved.\medskip

$(4)\Rightarrow(3)\;$ For any two distinct points $a$ and $d\in Z$, let $c$ be a fixed point in $(Z, \rho)$, which is different from $a$ and $d$. Then there is a point $x_0\in Z$ such that $$\mu_1\leq r(a,x_0,c,d)\leq \mu_2,$$
where $0< \mu_1\leq \mu_2<1$.

 By repeating this procedure, we can find a sequence $\{x_i\}_{i\in \mathbb{N}^+}$ in $(Z, \rho)$ such that $$\mu_1\leq r(a,x_i,x_{i-1},d)\leq \mu_2.$$ Then $$\mu_1\leq \frac{r(a,x_i,c,d)}{r(a,x_{i-1},c,d)}=r(a,x_i,x_{i-1},d)\leq\mu_2,$$ which implies that $$\mu_1^{i+1}\leq r(a,x_i,c,d)\leq \mu_2^{i+1},$$ and so $x_i\rightarrow d$ as $i\rightarrow+\infty$.

Similarly, we know that there exists $\{x_{-i}\}_{i\in \mathbb{N}_+}$ in $(Z, \rho)$ such that $$\mu_1\leq r(d,x_{-i},x_{1-i},a)=r(a,x_{1-i},x_{-i},d)\leq \mu_2$$
and
$$\mu_1\mu_2^{1-i}\leq r(a,x_{-i},c,d)\leq \mu_1^{1-i}\mu_2.$$
Then $x_{-i}\rightarrow a$ as $i\rightarrow+\infty,$ and hence we have proved that $(Z, \rho)$ is $\frac{1}{\mu_1}$-dense.
\qed

%%%%%%%%%%%%%%%%%%%%%%%%%%%%%%%
%%%%%%%%%%%%%%%%%%%%%%%%%%%%%%%
\section{The invariant property of uniform perfectness with respect to quasim\"obius maps}\label{sec-4}
%%%%%%%%%%%%%%%%%%%%%%%%%%%%%%%
%%%%%%%%%%%%%%%%%%%%%%%%%%%%%%%

The aim of this section is to prove Theorem \ref{lem3}. To this end, by Theorem \ref{main-thm1}, it suffices to show the following lemma.

\begin{lem}\label{lem-3} Let $f:$ $(Z_1,\rho_1)\to (Z_2,\rho_2)$ be $\theta$-quasim\"obius, where both $(Z_i,\rho_i)$ ($i=1$, $2$) are quasi-metric. Then $(Z_1, \rho_1)$ is $\sigma$-dense if and only if $(Z_2,\rho_2)$ is $\sigma'$-dense, quantitatively.
\end{lem}
\bpf Since the inverse of a $\theta$-quasim\"obius map is $\theta'$-quasim\"obius with $\theta'(t)=\frac{1}{\theta^{-1}(1/t)}$, we see that, to prove this lemma, it suffices to show that $(Z_2,\rho_2)$ is  $\sigma'$-dense under the assumption ``$(Z_1, \rho_1)$ being $\sigma$-dense", where $\sigma>1$ and $\sigma'$ depends only on $\sigma$ and $\theta$. For this, we only need to check that for each pair of points $a',$ $b'$ in $(Z_2, \rho_2)$, there is a $\sigma'$-chain in $(Z_2, \rho_2)$ joining them. Now, we assume that $\{x_i\}_{i\in \mathbb{Z}}$ is a $\sigma$-chain in $(Z_1, \rho_1)$ joining the points $a$ and $b$
with $$\frac{1}{\sigma}\leq r(a,x_i,x_{i+1},b)\leq \sigma.$$ Then for all $i$, we have $$\frac{1}{\theta(\sigma)+1}\leq r(a',x_i',x_{i+1}',b')\leq \theta(\sigma)+1,$$ which shows that $\{x_i'\}_{i\in \mathbb{Z}}$ is a $\sigma'$-chain in $(Z_2, \rho_2)$ joining $a'$ and $b'$ with $\sigma'=\theta(\sigma)+1$.
\epf

%%%%%%%%%%%%%%%%%%%%%%%%%%%%%%%
%%%%%%%%%%%%%%%%%%%%%%%%%%%%%%%
\section{uniform perfectness, (power) quasisymmetric maps and (power) quasim\"obius maps}\label{sec-5}
%%%%%%%%%%%%%%%%%%%%%%%%%%%%%%%
%%%%%%%%%%%%%%%%%%%%%%%%%%%%%%%

This section is devoted to the proof of Theorem \ref{main-thm0} concerning the relationships among uniform perfectness, (power) quasisymmetric maps and (power) quasim\"obius maps in quasi-metric spaces. It consists of two subsections. In the first subsection, we shall prove a relationship among uniform perfectness, quasisymmetric maps and power quasisymmetric maps, i.e. Theorem \ref{main-thm0}\eqref{main-1}, and in the second subsection, the proof of a relationship among uniform perfectness, quasim\"obius maps and power quasim\"obius maps, i.e. Theorem \ref{main-thm0}\eqref{main-2}, will be presented.

%%%%%%%%%%%%%%%%%%%%%%%%%%%%%%%
\subsection{The proof of Theorem \ref{main-thm0}\eqref{main-1}}
%%%%%%%%%%%%%%%%%%%%%%%%%%%%%%%

Let $\varepsilon\in (0,1)$ be a constant such that $K^{\varepsilon}\leq 2$. Then it follows from Lemma \ref{lem4} that
 there exists a metric $d_{\varepsilon}$ (briefly $d$ in the following) in $Z$ such that $$\frac{1}{4}\rho^{\varepsilon}(z_1,z_2)\leq d(z_1,z_2)\leq \rho^{\varepsilon}(z_1,z_2) $$ for all $z_1,$ $z_2\in Z$. Let $id$ denote the identity map from $(Z,\rho)$ to $(Z,d)$, i.e.,
 $$id:\; (Z, \rho)\to (Z, d).$$ Obviously, $id$ is power quasisymmetric with its control function $\eta(t)=4(t^{\epsilon}\vee t^\frac{1}{\epsilon})$.

  We first assume that $(Z,\rho)$ is uniformly perfect, and consider a quasisymmetric map $f$ defined in $(Z,\rho)$.
 It follows from the power quasisymmetry of $id$ and Lemma \ref{lem2} that $(Z,d)$ is uniformly perfect, and so Theorem \ref{main-thm1} implies that $(Z,d)$ is $(\lambda_1, \lambda_2)$-HD for constants $\lambda_1$ and $\lambda_2$ with $0<\lambda_1\leq \lambda_2<1$. Since $f\circ id^{-1}$ is quasisymmetric in $(Z,d)$, we see from \cite[Corollary 3.12]{TV} that $f\circ id^{-1}$ is power quasisymmetric, which implies that $f$ itself is power quasisymmetric.

Next, we assume that every quasisymmetric map of $(Z,\rho)$ is power quasisymmetric. Then we see that for any quasisymmetric map $g$ in $(Z,d)$, $g\circ id$ is power quasisymmetric in $(Z,\rho)$, and so $g$ itself is power quasisymmetric. Hence, by \cite[Theorems 4.13 and 6.20]{TroVai}, $(Z,d)$ is uniformly perfect. Since $id$ is power quasisymmetric, it follows from Lemma \ref{lem2} that $(Z, \rho)$ is uniformly perfect.
\qed

%%%%%%%%%%%%%%%%%%%%%%%%%%%%%%%
\subsection{The proof of Theorem \ref{main-thm0}\eqref{main-2}}
%%%%%%%%%%%%%%%%%%%%%%%%%%%%%%%

We start this subsection with the following result in metric spaces.

\begin{lem}\label{lem5}
Suppose $(Z,d)$ is a metric space with no isolated points. Then the following statements are quantitatively equivalent.
\ben
\item\label{wen-1}
$(Z,d)$ is uniformly perfect;
\item\label{wen-2}
every quasim\"obius map of $(Z,d)$ is power quasim\"obius.
\een
\end{lem}
\bpf
By \cite[Theorem 3.2]{Ase}, we only need to prove the implication from \eqref{wen-2} to \eqref{wen-1}. Assume that every quasim\"obius map in $(Z, d)$ is a power quasim\"obius map. To prove the uniform perfectness of $(Z, d)$, we divide the proof into two cases.
\bca\label{case1} $(Z,d)$ is unbounded.\eca

We shall apply Theorem \ref{main-thm0}\eqref{main-1} to finish the proof in this case. For this, we assume that $f$ is a quasisymmetric map in $(Z,d)$. Then Lemma \ref{w-02} implies that $f$ is quasim\"obius, and further, \cite[Theorem 3.10]{Vai89} guarantees that $f(z)\to \infty$ as $z\to \infty$. Again, it follows from \cite[Theorem 3.10]{Vai89} that $f$ is power quasisymmetric, and so Theorem \ref{main-thm0}\eqref{main-1} ensures that $(Z,d)$ is uniformly perfect. Hence the lemma is true in this case.
\bca\label{case2}$(Z,d)$ is bounded.\eca

By the Kuratowski embedding theorem \cite{Kur}, we may assume that $Z$ is a subset of a Banach space $E$.
By performing an auxiliary translation, further, we assume that $0\in Z$. Let  $$u(x)=\frac{x}{|x|^2}$$
 be the inversion in $\dot{E}=E\cup\{\infty\}$. Then, clearly, $u(Z)$ is unbounded. By \cite[\S 1.6]{Vai89}, $u$ is $\theta$-quasim\"obius, where $\theta(t)=81t$, and obviously, it is power quasim\"obius. To prove that $(Z, d)$ is uniformly perfect, by Theorem \ref{lem3}, it suffices to show that $u(Z)$ is uniformly perfect. Again, we shall apply Theorem \ref{main-thm0}\eqref{main-1} to reach this goal. For this, we assume that $g$ is quasisymmetric in $u(Z)$. Once more, by Lemma \ref{w-02}, $g$ is quasim\"obius. Then $g\circ u$ is quasim\"obius in $(Z,d)$, which implies that $g\circ u$ is power quasim\"obius, and thus we deduce from Lemma \ref{sun-4}\eqref{sun-6} that $g$ itself is power quasim\"obius. So we infer from \cite[Theorem 3.10]{Vai89} that $g$ is power quasisymmetric. Then it follows from Theorem \ref{main-thm0}\eqref{main-1} that $u(Z)$ is uniformly perfect. Hence the proof of the lemma is complete.
\epf

\noindent {\bf The proof of Theorem \ref{main-thm0}\eqref{main-2}}.\quad Let $id:$ $(Z, \rho)\to (Z, d)$ be the same as that in the proof of Theorem \ref{main-thm0}\eqref{main-1}. Then $id$ is power quasisymmetric with its control function $\eta(t)=4(t^{\epsilon}\vee t^\frac{1}{\epsilon})$, where $\epsilon\in (0,1)$.

Assume now that $(Z, \rho)$ is uniformly perfect, and so is $(Z, d)$ by Lemma \ref{lem2}. For any quasim\"obius map $f$ in $(Z, \rho)$, it follows from Lemma \ref{lem5} that $f\circ id$ is power quasim\"obius, and so is $f$ itself by Lemma \ref{sun-4}. This shows that the necessity in Theorem \ref{main-thm0}\eqref{main-2} is true.

To prove the sufficiency in Theorem \ref{main-thm0}\eqref{main-2}, it suffices to prove the uniform perfectness of $(Z, d)$ under the assumption that every quasisymmetric map in $(Z, \rho)$ is power quasim\"obius. By Lemma \ref{lem5}, we only need to show the power quasisymmetry of each quasisymmetric map in $(Z, d)$. This fact easily follows from Lemma \ref{sun-4}. Hence the proof of Theorem \ref{main-thm0}\eqref{main-2} is complete.\qed

\section{Applications}\label{sec-6}

The aim of this section is twofold. First, as an application of Theorem \ref{lem3}, we will give a different proof to \cite[Theorem 7.1]{Me}. Second, we shall apply Theorem \ref{lem3} to discuss the uniform perfectness of a complete quasi-metric space and the corresponding boundary of its hyperbolic approximation.

\subsection{Application $I$}
We begin this subsection with a definition.

\bdefe
For $p\in (Z,\rho)$, let
$$\rho_p(x,y)=\frac{r^2\rho(x,y)}{\rho(x,p)\rho(y,p)}$$
for all $x,$ $y\in Z\setminus\{p\}$. Then $\rho_p$ is said to be {\it the inversion} with respect to $\rho$ centered at $p$ with radius $r>0$.
\edefe

\begin{thm}\label{app-1} $($\cite[Theorem 7.1]{Me}$)\;$
 For any $p\in Z$, if
 $(Z\setminus\{p\},\rho)$ is a uniformly perfect quasi-metric space, then $(Z\setminus\{p\},\rho_p)$ is a uniformly perfect quasi-metric space.\end{thm}

\bpf First, if $(Z\setminus\{p\},\rho)$ is a $K$-quasi-metric space, by \cite[Proposition 5.3.6]{BuSc}, we know that $(Z\setminus\{p\},\rho_p)$ is a $K^2$-quasi-metric space. Then a direct computation gives that the identity map from $(Z\setminus\{p\},\rho)$ to $(Z\setminus\{p\},\rho_p)$ is $\theta$-quasim\"{o}bius with $\theta(t)=t$. Hence the proof of the theorem easily follows from Theorem \ref{lem3}.
\epf

\subsection{Application $II$}

Let $Hyp_r(Z,\rho)$ denote the hyperbolic approximation of $(Z, \rho)$ with parameter $r$,  $\partial_{\infty}^{a,o}Hyp_r(Z,\rho)$ the boundary at infinity of $Hyp_r(Z,\rho)$ with respect to the quasi-metric $a^{-(\cdot|\cdot)_o}$ based at $o\in Hyp_r(Z,\rho)$ with $a>1$, and $\partial_{\infty}^{a',b}Hyp_r(Z,\rho)$ the boundary at infinity of $Hyp_r(Z,\rho)$ with respect to the quasi-metric $a'^{-(\cdot|\cdot)_\omega}$ based at $\omega$ with $a'>1$, where $b$ is a Busemann function based at $\omega$. See \cite[\S$3$]{JJ} for their precise definitions.

\begin{thm}\label{app-2}
Suppose $(Z,\rho)$ is a complete quasi-metric space and $r\in (0,1)$.
Then the following are quantitatively equivalent.\ben
\item[$(a)$]
$(Z,\rho)$ is uniformly perfect;
\item[$(b)$]
$\partial_{\infty}^{a,o}Hyp_r(Z,\rho)$ is uniformly perfect;
\item[$(c)$]
$\partial_{\infty}^{a',b}Hyp_r(Z,\rho)$ is uniformly perfect. \een\end{thm}

\bpf First, by \cite[Proposition 2.2.9 and 5.2.8]{BuSc}, we know that the identity map from $\partial_{\infty}^{a',b}Hyp_r(Z,\rho)$ to $\partial_{\infty}^{a,o}Hyp_r(Z,\rho)$ is quasim\"obius, and so Theorem \ref{lem3} implies the quantitative equivalence of $(b)$ and $(c)$.

To finish the proof of the theorem, we divide the discussions into two cases. The first case is that $(Z,\rho)$ is unbounded. By \cite[Theorem 3]{JJ}, we know that for any Busemann function $b\in B(\omega)$, the identity map from $\partial_{\infty}^{a',b}Hyp_r(Z,\rho)$ to $(Z,\rho)$ is
bi-H\"older, and thus Theorem \ref{lem3} guarantees the quantitative equivalence of $(a)$ and $(c)$. For the remainder case, that is, $(Z,\rho)$ is bounded, again, by \cite[Theorem 3]{JJ}, we see that the identity map from $\partial_{\infty}^{a,o}Hyp_r(Z,\rho)$ to $(Z,\rho)$ is  bi-H\"older. Once more, it follows from Theorem \ref{lem3} that $(a)$ and $(b)$ are quantitatively equivalent. Hence the proof of this theorem is complete.
\epf

%{\bf Acknowledgement.}
\bigskip
{\bf Acknowledgements.}
The authors would like to thank Professors Xiantao Wang and  Manzi Huang for several comments on this manuscripts.

\end{document}